\documentclass[12pt,reqno]{amsart}

\newcommand\version{May 3, 2012}

\usepackage{amsmath,amsfonts,amsthm,amssymb,amsxtra}

\newtheorem{theorem}{Theorem}[section]
\newtheorem{proposition}[theorem]{Proposition}
\newtheorem{lemma}[theorem]{Lemma}
\newtheorem{corollary}[theorem]{Corollary}

\theoremstyle{definition}

\theoremstyle{remark}

\numberwithin{equation}{section}

\newcommand{\C}{\mathbb{C}}

\renewcommand{\epsilon}{\varepsilon}

\renewcommand{\phi}{\varphi}
\newcommand{\R}{\mathbb{R}}

\newcommand{\Sph}{\mathbb{S}}

\begin{document}

\title[Sharp HLS inequality --- \version]{A new, rearrangement-free proof of the \\ sharp Hardy-Littlewood-Sobolev inequality}

\dedicatory{Dedicated to D. E. Edmunds and W. D. Evans}

\author{Rupert L. Frank}
\address{Rupert L. Frank, Department of Mathematics, Princeton University, Princeton, NJ 08544, USA}
\email{rlfrank@math.princeton.edu}

\author{Elliott H. Lieb}
\address{Elliott H. Lieb, Departments of Mathematics and Physics, Princeton University, P.~O.~Box 708, Princeton, NJ 08544, USA}
\email{lieb@math.princeton.edu}

\thanks{\copyright\, 2012 by the authors. This paper may be reproduced, in its entirety, for non-commercial purposes.\\
Support by U.S. NSF grant PHY 0965859 (E.H.L.) is acknowledged.}

\begin{abstract}
 We show that the sharp constant in the Hardy-Littlewood-Sobolev inequality can be derived using the method that we employed earlier for a similar inequality on the Heisenberg group. The merit of this proof is that it does not rely on rearrangement inequalities; it is the first one to do so for the whole parameter range.
\end{abstract}

\maketitle

\section{Introduction}

In a recent paper \cite{FrLi2} we showed how to compute the sharp
constants for the analogue of the Hardy-Littlewood-Sobolev (HLS)
inequality on the Heisenberg group. Unlike the situation for the usual
HLS inequality on $\R^N$, there is no known useful symmetric
decreasing rearrangement technique for the Heisenberg group analogue.
A radically new approach had to be developed and that approach can, of
course, be used for the original HLS problem as well, thereby
providing a genuinely rearrangement-free proof of HLS on $\R^N$. That
will be given here.

The HLS inequality (more precisely, the diagonal case) on $\R^N$ is
\begin{equation} \label{eq:hls}
\boxed{ \phantom{\Big|} \Big| \iint_{\R^N\times\R^N}\! \frac{\overline{f(x)}\ g(y)}{|x-y|^\lambda} \,dx\,dy \Big|
\leq 
\pi^{\lambda/2} \frac{\Gamma((N-\lambda)/2)}
{\Gamma(N-\lambda/2)} \left( \frac{\Gamma(N)}{\Gamma(N/2)} \right)^{1-\lambda/N} \! \|f\|_p \ \|g\|_p }
\end{equation}
where $0<\lambda<N$ and $p=2N/(2N-\lambda)$. The constant in \eqref{eq:hls} is sharp
and inequality \eqref{eq:hls} is strict unless $f$ and $g$ are
proportional to a common translate or dilate of
\begin{equation}
 \label{eq:opt}
H(x)= \left(1+|x|^2\right)^{-(2N-\lambda)/2}\,.
\end{equation}

An equivalent formulation of \eqref{eq:hls}, which has been noted before in the special cases $s=1$ and $s=1/2$ \cite[Thms. 8.3 and 8.4]{LiLo}, is the sharp fractional Sobolev inequality\footnote{The published version of this paper contains a typo in the following boxed formula \eqref{eq:fracsob}, which has been corrected here. We thank T. Weth for pointing this out to us.}
\begin{equation}
 \label{eq:fracsob}
\boxed{
\ \left\| (-\Delta)^{s/2} u \right\|^2 \geq \frac{2^{2s} \pi^s \, \Gamma((N+2s)/2)}{\Gamma((N-2s)/2)} \left(\frac{\Gamma(N/2)}{\Gamma(N)}\right)^{2s/N} \|u\|_q^2 \ }
\end{equation}
for $0<s<N/2$ and $q=2N/(N-2s)$. This follows from \eqref{eq:hls} by a duality argument (see \cite[Thm. 8.3]{LiLo}), using the fact that the Green's function of $(-\Delta)^s$ is $2^{-2s}\pi^{-N/2} \Gamma((N-2s)/2)/(\Gamma(s)|x|^{N-2s})$ for $0<s<N/2$ \cite[Thm. 5.9]{LiLo}. In particular, for $s=1$, \eqref{eq:fracsob} is the familiar Sobolev inequality
\begin{equation}
 \label{eq:sobjl}
\int_{\R^N} |\nabla u |^2 \,dx
\geq \frac{N(N-2)}{4} \left(\frac{2\pi^{(N+1)/2}}{\Gamma((N+1)/2)}\right)^{2/N} \ \|u\|_q^2
\end{equation}
for $N\geq 3$ and $q=2N/(N-2)$ in the sharp form of \cite{Ro,Au,Ta}.

To recall, briefly, the previous proofs of \eqref{eq:hls} we first mention
the papers \cite{HaLi1,HaLi2,So}, where the inequality was initially derived,
but with a non-sharp constant. The sharp version was found in \cite{Li} by
noting the conformal invariance of the problem and relating it, via
stereographic projection, to a conformally equivalent, but more tractable
problem on the sphere $\Sph^N$. Riesz's rearrangement inequality (see
\cite[Thm. 3.7]{LiLo}) was used in the proof of the existence of a
maximizer, and its strong version (\cite{Li0}, see also \cite[Thm. 3.9]{LiLo}) was used to prove
that the constant function is a maximizer -- in the spherical version.
There are other, by now standard,  ways to prove the existence of  a 
maximizer; that is not the issue. The main point is to prove that
(\ref{eq:opt}) is a maximizer and that it is, essentially, unique.

Then Carlen and Loss \cite{CaLo} cleverly utilized the translational
symmetry in $\R^N$ in competition with the rotational symmetry on the
sphere, together with the strong Riesz inequality, to conclude the same
thing.

Another proof, but only for $N-2\leq\lambda < N$, was recently given in
\cite{FrLi1}. This was done by proving a form of reflection positivity
for inversions in spheres in $\R^N$, and generalizing a theorem of
Li and Zhu \cite{LZh}.  This is the first rearrangement-free proof of HLS, but it is
not valid for $0<\lambda < N-2$. An elegant, rearrangement-free proof, this time only for $\lambda=N-2$, is in \cite{CCaLo}.

In this note we show how the new method developed in \cite{FrLi2} can be
adapted to the HLS problem to yield a proof for all $0<\lambda <N$. We also apply the method directly to a proof of \eqref{eq:sobjl} in Section \ref{sec:jl}.


\section{Main result}

We shall prove

\begin{theorem}\label{main}
 Let $0<\lambda<N$ and $p:=2N/(2N-\lambda)$. Then \eqref{eq:hls} holds for any $f,g\in L^p(\R^N)$.
Equality holds if and only if
$$
f(x) = c \ H(\delta(x-a)) \,,
\qquad
g(y) = c' \ H(\delta(x-a))
$$
for some $c,c'\in\C$, $\delta>0$ and $a\in\R^N$ (unless $f\equiv 0$ or $g\equiv 0$). Here $H$ is the function in \eqref{eq:opt}.
\end{theorem}

In other words, we prove that the function $H$ in \eqref{eq:opt} is the unique optimizer in inequality \eqref{eq:hls} up to translations, dilations and multiplication by a constant.

The stereographic projection (see Appendix \ref{sec:equiv}) defines a bijection between $\R^N$ and the punctured sphere $\Sph^N\setminus\{(0,\ldots,0,-1)\}$. We consider the sphere $\Sph^N$ as a subset of $\R^{N+1}$ with coordinates $(\omega_1,\ldots,\omega_{N+1})$ satisfying $\sum_{j=1}^{N+1} \omega_j^2 =1$, and (non-normalized) measure denoted by $d\omega$. Via stereographic projection Theorem \ref{main} is equivalent to

\begin{theorem}\label{mainsph}
 Let $0<\lambda<N$ and $p:=2N/(2N-\lambda)$. Then for any $f,g\in L^p(\Sph^N)$
\begin{equation}
 \label{eq:mainsph}
\left| \iint_{\Sph^N\times\Sph^N} \frac{\overline{f(\omega)}\ g(\eta)}{|\omega-\eta|^{\lambda}} \,d\omega\,d\eta \right|
\leq \pi^{\lambda/2} \frac{\Gamma((N-\lambda)/2)}{\Gamma(N-\lambda/2)} \left( \frac{\Gamma(N)}{\Gamma(N/2)} \right)^{1-\lambda/N}
 \|f\|_p\ \|g\|_p
\end{equation}
with equality if and only if
\begin{equation}
 \label{eq:optsph}
f(\omega) = \frac{c}{(1-\xi\cdot\omega)^{(2N-\lambda)/2}} \,,
\qquad
g(\omega) = \frac{c'}{(1-\xi\cdot\omega)^{(2N-\lambda)/2}} \,,
\end{equation}
for some $c,c'\in\C$ and some $\xi\in\R^{N+1}$ with $|\xi|<1$ (unless $f\equiv 0$ or $g\equiv 0$).
\end{theorem}

In particular, with $\xi=0$, $f=g\equiv1$ are optimizers.

We conclude this section by recalling that \eqref{eq:mainsph} can be differentiated at the endpoints $\lambda=0$ and $\lambda=N$, where the inequality turns into an equality. In this way one obtains the logarithmic HLS inequality \cite{CaLo2,Be2} and a conformally invariant logarithmic Sobolev inequality \cite{Be1}.


\section{The sharp Sobolev inequality on the sphere}\label{sec:jl}

In this section we derive the classical Sobolev inequality \eqref{eq:sobjl}. This case is simpler than the general $\lambda$ case of the HLS inequality, but it already contains the main elements of our strategy. It is easiest for us to work in the formulation on the sphere $\Sph^{N}$.

We consider $\Sph^N$ as a subset of $\R^{N+1}$, i.e., $\{(\omega_1,\ldots,\omega_{N+1}):\ \sum_{j=1}^{N+1} \omega_j^2 =1\}$. We recall that the conformal Laplacian on $\Sph^N$ is defined by
$$
\mathcal L := -\Delta + \frac{N(N-2)}4 \,,
$$
where $\Delta$ is the Laplace-Beltrami operator on $\Sph^{N}$, and we denote the associated quadratic form by
$$
\mathcal E[u] := \int_{\Sph^N} \left( |\nabla u|^2 + \frac{N(N-2)}4 |u|^2 \right) \,d\omega \,.
$$
The sharp Sobolev inequality on $\Sph^N$ is

\begin{theorem}\label{mainjl}
For all $u\in H^1(\Sph^N)$ one has
\begin{equation}\label{eq:mainjl}
\mathcal E[u] \geq \frac{N(N-2)}{4} \left(\frac{2\pi^{(N+1)/2}}{\Gamma((N+1)/2)}\right)^{2/N} \left( \int_{\Sph^N} |u|^{2N/(N-2)} \,d\omega \right)^{(N-2)/N} \,,
\end{equation}
with equality if and only if
\begin{equation}
 \label{eq:optjl}
u(\omega) = c \, (1-\xi\cdot\omega)^{-(N-2)/2}
\end{equation}
for some $c\in\C$ and some $\xi\in\R^{N+1}$ with $|\xi|<1$.
\end{theorem}

See Appendix \ref{sec:equiv} for the equivalence of the $\R^N$-version \eqref{eq:sobjl} and the $\Sph^N$-version \eqref{eq:mainjl} of the Sobolev inequality.

In the proof of Theorem \ref{mainjl} we shall make use of the following elementary formula.

\begin{lemma}\label{gsr}
 For all $u\in H^1(\Sph^N)$ one has
\begin{equation}
\label{eq:gsr}
\sum_{j=1}^{N+1} \mathcal E[\omega_j u] = \mathcal E[u] + N \int_{\Sph^N} |u|^2 \,d\omega \,.
\end{equation}
\end{lemma}

\begin{proof}
We begin by noting that for any smooth, real-valued function $\phi$ on $\Sph^N$ one has
\begin{align*}
|\nabla (\phi u)|^2 = \phi^2 |\nabla u|^2 + |u|^2 |\nabla\phi|^2 + \phi \nabla\phi\cdot \nabla(|u|^2) \,.
\end{align*}
Hence an integration by parts leads to
\begin{align*}
\int_{\Sph^N} |\nabla (\phi u)|^2 \,d\omega = \int_{\Sph^N} \left( \phi^2 |\nabla u|^2 - \phi(\Delta\phi) |u|^2 \right) \,d\omega  \,.
\end{align*}
We apply this identity to $\phi(\omega)=\omega_j$. Using the fact that
$$
-\Delta \omega_j = N \omega_j \,,
$$
we find
\begin{align*}
\int_{\Sph^N} |\nabla (\omega_j u)|^2 \,d\omega = \int_{\Sph^N} \omega_j^2 \left(  |\nabla u|^2 + N |u|^2 \right) \,d\omega  \,.
\end{align*}
Summing over $j$ yields \eqref{eq:gsr} and completes the proof.
\end{proof}

We are now ready to give a short

\begin{proof}[Proof of Theorem \ref{mainjl}]
It is well-known that there is an optimizer $U$ for inequality \eqref{eq:mainjl}. (Using the stereographic projection, one can deduce this for instance from the existence of an optimizer on $\R^N$; see \cite{Li}.)

As a preliminary remark we note that any optimizer is a complex multiple of a non-negative function. Indeed, if $u=a+ib$ with $a$ and $b$ real functions, then $\mathcal E[u]=\mathcal E[a]+\mathcal E[b]$. We also note that the right side of \eqref{eq:mainjl} is $\|a^2+b^2\|_{q/2}$ with $q=2N/(N-2)>2$. By the triangle inequality, $\|a^2+b^2\|_{q/2}\leq \|a^2\|_{q/2}+\|b^2\|_{q/2}$. This inequality is strict unless $a\equiv 0$ or $b^2=\lambda^2 a^2$ for some $\lambda\geq 0$. Therefore, if $U=A+iB$ is an optimizer for \eqref{eq:mainjl}, then either one of $A$ and $B$ is identically equal to zero or else both $A$ and $B$ are optimizers and $|B|=\lambda |A|$ for some $\lambda>0$. For any real $u\in H^1(\Sph^N)$ its positive and negative parts $u_\pm$ belong to $H^1(\Sph^N)$ and satisfy $\partial u_\pm/\partial\omega_k = \pm\chi_{\{\pm u>0\}} \partial u/\partial\omega_k$ in the sense of distributions. (This can be proved similarly to \cite[Thm. 6.17]{LiLo}.) Thus $\mathcal E[u]=\mathcal E[u_+]+\mathcal E[u_-]$ for real $u$. Moreover, $\| u\|_q^2 \leq \|u_+\|_q^2 + \|u_-\|_q^2$ for real $u$ with strict inequality unless $u$ has a definite sign. Therefore, if $U=A+iB$ is an optimizer for \eqref{eq:mainjl}, then both $A$ and $B$ have a definite sign. We conclude that any optimizer is a complex multiple of a non-negative function. Hence we may assume that $U\geq 0$.

It is important for us to know that we may confine our search for optimizers to functions $u$ satisfying the `center of mass condition'
\begin{equation}\label{eq:comjl}
\int_{\Sph^N} \omega_j\ |u(\omega)|^q \,d\omega= 0\,,
\qquad j=1,\ldots,N+1 \,.
\end{equation}
It is well-known, and used in many papers on this subject (e.g., \cite{He,On,ChYa}), that \eqref{eq:comjl} can be assumed, and we give a proof of this fact in Appendix \ref{sec:com}. It uses three facts: one is that inequality \eqref{eq:mainjl} is invariant under $O(N+1)$ rotations of $\Sph^N$. The second is that the stereographic projection, that maps $\R^N$ to $\Sph^N$, leaves the optimization problem invariant. The third is that the $\R^N$-version, \eqref{eq:sobjl}, of inequality \eqref{eq:mainjl} is invariant under dilations $F(x) \mapsto \delta^{(N-2)/2} F(\delta x)$. Our claim in the appendix is that by a suitable choice of $\delta$ and a rotation we can achieve \eqref{eq:comjl}.

Therefore we may assume that the optimizer $U$ satisfies \eqref{eq:comjl}. Imposing this constraint does not change the positivity of $U$. We shall prove that the only optimizer with this property is the constant function (which leads to the stated expression for the sharp constant). It follows, then, that the only optimizers without condition \eqref{eq:comjl} are those functions for which the dilation and rotation, just mentioned, yields a constant. In Appendix \ref{sec:com} we identify those functions as the functions stated in \eqref{eq:optjl}.

The second variation of the quotient $\mathcal E[u]/\|u\|_q^2$ around $u=U$ shows that
\begin{equation}
 \label{eq:secvarjl}
\mathcal E[v] \int_{\Sph^N} U^q \,d\omega - (q-1) \mathcal E[U] \int_{\Sph^N} U^{q-2} |v|^2 \,d\omega \geq 0
\end{equation}
for all $v$ with $\int U^{q-1} v \,d\omega =0$.

Because $U$ satisfies condition \eqref{eq:comjl} we may choose $v(\omega)=\omega_j U(\omega)$ in \eqref{eq:secvarjl} and sum over $j$. We find
\begin{equation}\label{eq:secondvarjl}
\sum_{j=1}^{N+1} \mathcal E[\omega_j U] \geq (q-1)\ \mathcal E[U] \,.
\end{equation}
On the other hand, Lemma \ref{gsr} with $u=U$ implies
\begin{equation*}
\sum_{j=1}^{N+1} \mathcal E[\omega_j U] = \mathcal E[U] + N \int_{\Sph^N} U^2 \,d\omega \,,
\end{equation*}
which, together with \eqref{eq:secondvarjl}, yields
$$
N \int_{\Sph^N} U^2 \,d\omega \geq (q-2)\ \mathcal E[U] \,.
$$
Recalling that $q-2=\frac{4}{N-2}$, we see that this is the same as
$$
\int_{\Sph^N} |\nabla U|^2 \,d\omega \leq 0 \,.
$$
We conclude that $U$ is the constant function, as we intended to prove.
\end{proof}



\section{The sharp HLS inequality on the sphere}\label{sec:proofs}

Our goal in this section is to compute the sharp constant in inequality \eqref{eq:mainsph} on the sphere $\Sph^N$. We outline our argument in Subsection \ref{sec:strategy} and reduce everything to the proof of a \emph{linear} inequality. After some preparations in Subsection \ref{sec:funk} we shall prove this inequality in Subsection~\ref{sec:mainproof}.


\subsection{Strategy of the proof}\label{sec:strategy}

\emph{Step 1.} The optimization problem corresponding to \eqref{eq:mainsph} admits an optimizing pair with $f=g$. The fact that one only needs to consider $f=g$ follows from the positive definiteness of the kernel $|x-y|^{-\lambda}$. The existence of an optimizer has been proved in \cite{Li} for the inequality \eqref{eq:hls} on $\R^N$ and follows, as explained in Appendix \ref{sec:equiv}, via stereographic projection for the inequality on the sphere; for a rearrangement-free proof, see \cite{Ln} and also the arguments in \cite{FrLi2}, which easily carry over to the $\R^N$ case.

We claim that any optimizer for problem \eqref{eq:mainsph} with $f=g$ is a complex multiple of a non-negative function. Indeed,
if we denote the left side of \eqref{eq:hls} with $g=f$ by $I[f]$ and if $f=a+ib$ for real functions $a$ and $b$, then $I[f]=I[a]+I[b]$. Moreover, for any numbers $\alpha,\beta,\gamma,\delta\in\R$ one has $\alpha\gamma+\beta\delta \leq \sqrt{\alpha^2+\beta^2}\sqrt{\gamma^2+\delta^2}$ with strict inequality unless $\alpha\gamma+\beta\delta\geq 0$ and $\alpha\delta=\beta\gamma$. Since the kernel $|x-y|^{-\lambda}$ is strictly positive, we infer that $I[a]+I[b]\leq I[\sqrt{a^2+b^2}]$ for any real functions $a,b$ with strict inequality unless $a(x)a(y)+b(x)b(y)\geq 0$ and $a(x)b(y)=a(y)b(x)$ for almost every $x,y\in\R^N$. From this one easily concludes that any optimizer is a complex multiple of a non-negative function.

We denote a non-negative optimizer for problem \eqref{eq:mainsph} by $h:=f=g$. Since $h$ satisfies the Euler-Lagrange equation
$$
\int_{\Sph^N} \frac{h(\eta)}{{|\omega -\eta|^{\lambda}}} \,d\eta = c\, h^{p-1}(\omega) \,,
$$
we see that $h$ is \emph{strictly} positive.

\emph{Step 2.} As in the proof of Theorem \ref{mainjl}, we may assume that the center of mass of $h^p$ vanishes, that is,
\begin{equation}
 \label{eq:com}
\int_{\Sph^N} \omega_j \, h(\omega)^p \,d\omega = 0
\qquad \text{for} \ j=1,\ldots,N+1 \,.
\end{equation}
We shall prove that the only non-negative optimizer satisfying \eqref{eq:com} is the constant function. Then, for exactly the same reason as in the proof of Theorem \ref{mainjl}, the only optimizers without condition \eqref{eq:com} are the ones stated in \eqref{eq:optsph}. We also note that, once we know that a constant is the optimizer, the expression for the sharp constant follows by a computation (see the $l=0$ case of Corollary \ref{ev} below).

\emph{Step 3.} The second variation around the optimizer $h$ shows that
\begin{equation}
 \label{eq:secvar}
\iint \frac{\overline{f(\omega)}\ f(\eta)}{|\omega-\eta|^{\lambda}} \,d\omega\,d\eta 
\ \int h^{p} \,d\omega
- (p-1) \iint \frac{h(\omega)\ h(\eta)}{|\omega-\eta|^{\lambda}} \,d\omega\,d\eta
\ \int h^{p-2} |f|^2 \,d\omega \leq 0
\end{equation}
for any $f$ satisfying $\int h^{p-1} f \,d\omega =0$. Note that the term $h^{p-2}$ causes no problems (despite the fact that $p<2$) since $h$ is strictly positive.

Because of \eqref{eq:com} the functions $f(\omega)=\omega_j h(\omega)$ satisfy the constraint $\int h^{p-1} f \,d\omega =0$. Inserting them in \eqref{eq:secvar} and summing over $j$ we find
\begin{equation}
 \label{eq:secvarineq}
\iint \frac{h(\omega)\ \omega\cdot\eta \ h(\eta)}{|\omega-\eta|^{\lambda}} \,d\omega\,d\eta - (p-1) \iint \frac{h(\omega)\ h(\eta)}{|\omega-\eta|^{\lambda}} \,d\omega\,d\eta \leq 0 \,.
\end{equation}

\emph{Step 4.} \emph{This is the crucial step!} The proof of Theorem \ref{mainsph} is completed by showing that for \emph{any} (not necessarily maximizing) $h$ the inequality opposite to \eqref{eq:secvarineq} holds and is indeed strict unless the function is constant. This is the statement of the following theorem with $\alpha=\lambda/2$, noting that $p-1 = \alpha/(N-\alpha)$.

\begin{proposition}\label{keyineq}
 Let $0<\alpha<N/2$. For any $f$ on $\Sph^N$ one has
\begin{equation}\label{eq:keyineq}
 \iint \frac{\overline{f(\omega)}\ \omega\cdot\eta \ f(\eta)}{|\omega-\eta|^{2\alpha}} \,d\omega\,d\eta 
\geq \frac{\alpha}{N-\alpha} \iint \frac{\overline{f(\omega)}\ f(\eta)}{|\omega-\eta|^{2\alpha}} \,d\omega\,d\eta
\end{equation}
with equality iff $f$ is constant.
\end{proposition}

This proposition will be proved in Subsection \ref{sec:mainproof}.


\subsection{The Funk-Hecke theorem}\label{sec:funk}

We decompose $L^2(\Sph^N)$ into its $O(N+1)$-irreducible components,
\begin{equation}\label{eq:decomp}
L^2(\Sph^N)= \bigoplus_{l\geq 0} \mathcal H_l \,.
\end{equation}
The space $\mathcal H_l$ is the space of restrictions to $\Sph^N$ of harmonic polynomials on $\R^{N+1}$ which are homogeneous of degree $l$.

It is well-known that integral operators on $\Sph^N$ whose kernels have the form $K(\omega\cdot\eta)$ are diagonal with respect to this decomposition and their eigenvalues can be computed explicitly. A proof of the following Funk-Hecke formula can be found, e.g., in \cite[Sec. 11.4]{EMOT}. It involves the Gegenbauer polynomials $C_l^{(\lambda)}$, see \cite[Chapter 22]{AbSt}.

\begin{proposition}\label{funk}
 Let $K\in L^1((-1,1),(1-t^2)^{(N-2)/2}dt)$. Then the operator on $\Sph^N$ with kernel $K(\omega\cdot\eta)$ is diagonal with respect to decomposition \eqref{eq:decomp}, and on the space $\mathcal H_l$ its unique eigenvalue is given by
\begin{equation}
 \label{eq:funk}
\kappa_{N,l} \int_{-1}^1 K(t) C_l^{(N-1)/2}(t) (1-t^2)^{(N-2)/2} \,dt \,,
\end{equation}
where
\begin{equation*}
 \kappa_{N,l}=
\begin{cases}
 2 & \text{if}\ N=1\,,\ l=0\,,\\
l & \text{if}\ N=1\,,\ l\geq 1\,,\\
(4\pi)^{(N-1)/2} \, \frac{l!\ \Gamma((N-1)/2)}{(l+N-2)!} & \text{if}\ N\geq 2 \,.
\end{cases}
\end{equation*}
\end{proposition}


This proposition allows us to compute the eigenvalues of the family of operators appearing in Proposition \ref{keyineq}.

\begin{corollary}\label{ev}
Let $-1<\alpha<N/2$. The eigenvalue of the operator with kernel $(1-\omega\cdot\eta)^{-\alpha}$ on the subspace $\mathcal H_l$ is
\begin{equation}
 \label{eq:ev1}
E_l = \kappa_N \, 2^{-\alpha} \, (-1)^l \,\frac{\Gamma(1-\alpha)\, \Gamma(N/2-\alpha)}{\Gamma(-l+1-\alpha)\, \Gamma(l+N-\alpha)} \,,
\end{equation}
where
$$
\kappa_N =
\begin{cases}
2 \pi^{1/2}
& \text{if}\ N=1 \,, \\
2^{2(N-1)} \pi^{(N-1)/2} \, \frac{\Gamma((N-1)/2)\, \Gamma(N/2)}{(N-2)!} 
& \text{if}\ N\geq 2 \,.
\end{cases}
$$
When $\alpha$ is a non-negative integer, formula \eqref{eq:ev1} is to be understood by taking limits with fixed $l$.
\end{corollary}

This result appears already (without proof) in \cite{Be2}.

\begin{proof}
By Proposition \ref{funk} we have to evaluate the integral \eqref{eq:funk} for the choice $K(t) = (1-t)^{-\alpha}$. Our assertion follows from the $\beta=(N-2)/2-\alpha$ case of the formula
\begin{align}\label{eq:gegen}
 & \int_{-1}^1 (1+t)^{(N-2)/2} (1-t)^\beta C_l^{(N-1)/2}(t)\,dt  \\ 
& \quad = (-1)^l \, \frac{2^{N/2+\beta}\, \Gamma(1+\beta)\, \Gamma(N/2)\, \Gamma(l+N-1)\, \Gamma(-N/2+2+\beta)}{l!\, \Gamma(N-1)\, \Gamma(-l-N/2+2+\beta)\, \Gamma(l+N/2+1+\beta)} \,. \notag
\end{align}
This formula, which is valid for $\beta>-1$, follows from \cite[(7.311.3)]{GrRy} together with the fact that $C_l^{(\lambda)}(-t)=(-1)^l \, C_l^{(\lambda)}(t)$. As it stands, \eqref{eq:gegen} is only valid for $N\geq 2$. For $N=1$ and $l=0$, the (divergent) factors $\Gamma(l+N-1)$ and $\Gamma(N-1)$ need to be omitted, and for $N=1$ and $l\geq 1$, the divergent factor $\Gamma(N-1)$ in the denominator needs to be replaced by $\frac12$.
\end{proof}


\subsection{Proof of Proposition \ref{keyineq}}\label{sec:mainproof}

Using the fact that $|\omega-\eta|^2 = 2(1-\omega\cdot\eta)$, we see that the assertion is equivalent to
\begin{align*}
 & \iint \frac{\overline{f(\omega)}\,f(\eta)}{(1-\omega\cdot\eta)^{\alpha-1}} \,d\omega\,d\eta 
\leq \frac{N-2\alpha}{N-\alpha} \iint \frac{\overline{f(\omega)}\ f(\eta)}{(1-\omega\cdot\eta)^{\alpha}} \,d\omega\,d\eta \,.
\end{align*}
Both quadratic forms are diagonal with respect to decomposition \eqref{eq:decomp} and their eigenvalues on the subspace $\mathcal H_l$ are given by Corollary \ref{ev}. For simplicity, we first assume that $\alpha\neq 1$. The eigenvalue of the right side is $(N-2\alpha)E_l/(N-\alpha)$, with $E_l$ given by \eqref{eq:ev1}, and the eigenvalue of the left side is $\tilde E_l$, which is $E_l$ with $\alpha$ replaced by $\alpha-1$. Noting that
$$
\tilde E_l = E_l \frac{(\alpha-1) (N-2\alpha)}{(l-1+\alpha)(l+N-\alpha)}
$$
and that $E_l>0$ and $\alpha<N/2$, we see that the conclusion of the theorem is equivalent to the inequality
\begin{align*}
\frac{\alpha-1}{(l-1+\alpha)(l+N-\alpha)} \leq \frac{1}{N-\alpha}
\end{align*}
for all $l\geq 0$. This inequality is elementary to prove, distinguishing the cases $\alpha>1$ and $\alpha<1$. Finally, the case $\alpha=1$ is proved by letting $\alpha\to 1$ for fixed $l$.

Strictness of inequality \eqref{eq:keyineq} for non-constant $f$ follows from the fact that the above inequalities are strict unless $l=0$. This completes the proof of Proposition~\ref{keyineq}.
\qed



\appendix

\section{Equivalence of Theorems \ref{main} and \ref{mainsph}}\label{sec:equiv}

In this appendix we consider the stereographic projection $\mathcal S: \R^N \to \Sph^N$ and its inverse $\mathcal S^{-1}:  \Sph^{N}\to \R^N$ given by
\begin{align*}
 \mathcal S(x) = \left(\frac{2x}{1+|x|^2}, \frac{1-|x|^2}{1+|x|^2} \right) \,, \qquad
\mathcal S^{-1}(\omega) = \left(\frac{\omega_1}{1+\omega_{N+1}},\ldots,\frac{\omega_N}{1+\omega_{N+1}} \right) \,.
\end{align*}
The Jacobian of this transformation (see, e.g., \cite[Thm. 4.4]{LiLo}) is
$$
J_{\mathcal S}(x) = \left( \frac{2}{1+|x|^2} \right)^{N} \,,
$$
which implies that
\begin{equation}
 \label{eq:change}
\int_{\Sph^N} \phi(\omega) \,d\omega = \int_{\R^N} \phi(\mathcal S(x)) J_{\mathcal S}(x) \,dx
\end{equation}
for any integrable function $\phi$ on $\Sph^N$.

We now explain the equivalence of \eqref{eq:hls} and \eqref{eq:mainsph} for each fixed pair of parameters $\lambda$ and $p$ with $p=2N/(2N-\lambda)$. There is a one-to-one correspondence between functions $f$ on $\Sph^N$ and functions $F$ on $\R^N$ given by
\begin{equation}
 \label{eq:corres}
F(x)=|J_{\mathcal S}(x)|^{1/p} f(\mathcal S(x)) \,.
\end{equation}
It follows immediately from \eqref{eq:change} that $f\in L^p(\Sph^N)$ if and only if $F\in L^p(\R^N)$, and in this case $\|f\|_p = \|F\|_p$. Moreover, we note the fact that
$$
|\omega-\eta|^2 = \left( \frac{2}{1+|x|^2} \right) |x-y|^2 \left( \frac{2}{1+|y|^2} \right)
$$
for $\omega=\mathcal S(x)$ and $\eta=\mathcal S(y)$, where $|\omega-\eta|$ is the chordal distance between $\omega$ and $\eta$, i.e., the Euclidean distance in $\R^{N+1}$. With the help of this relation one easily verifies that
$$
\iint_{\R^N\times\R^N} \frac{\overline{F(x)}\ F(y)}{|x-y|^\lambda} \,dx\,dy 
= \iint_{\Sph^N\times\Sph^N} \frac{\overline{f(\omega)}\ f(\eta)}{|\omega-\eta|^{\lambda}} \,d\omega\,d\eta \,.
$$
This shows that the sharp constants in \eqref{eq:hls} and \eqref{eq:mainsph} coincide and that there is a one-to-one correspondence between optimizers. In particular, the function $f\equiv 1$ on $\Sph^N$ corresponds to the function
$$
|J_{\mathcal S}(x)|^{1/p} = 2^{N/p} H(x)
$$
on $\R^N$ with $H$ given in \eqref{eq:opt}.

Similarly, when $p=2N/(N-2)$, and $F$ and $f$ are related via \eqref{eq:corres}, then
\begin{equation}
 \label{eq:corressob}
\int_{\R^N} | \nabla F |^2 \,dx 
= \int_{\Sph^N} \left( |\nabla f|^2 +  \frac{N(N-2)}4 |f|^2 \right) \,d\omega \,,
\end{equation}
as can be checked by a direct computation.


\section{The center of mass condition}\label{sec:com}

Here, we prove that by a suitable inequality preserving transformation of $\Sph^N$ we may assume the center of mass conditions given in \eqref{eq:comjl} and \eqref{eq:com}.

We shall define a family of maps $\gamma_{\delta,\xi}:\Sph^N\to\Sph^N$ depending on two parameters $\delta>0$ and $\xi\in\Sph^N$. To do so, we denote dilation on $\R^N$ by $\mathcal D_\delta$, that is, $\mathcal D_\delta(x)=\delta x$. Moreover, for any $\xi\in\Sph^N$ we choose an orthogonal $(N+1)\times(N+1)$ matrix $O$ such that $O\xi=(0,\ldots,0,1)$ and we put
$$
\gamma_{\delta,\xi}(\omega) := O^T \mathcal S\left(\mathcal D_\delta\left(\mathcal S^{-1}\left(O\omega\right)\right)\right)
$$
for all $\omega\in\Sph^N\setminus\{-\xi\}$ and $\gamma_{\delta,\xi}(-\xi):=-\xi$. This transformation depends only on $\xi$ (and $\delta$) and not on the particular choice of $O$. Indeed, a straightforward computation shows that
\begin{equation*}
\gamma_{\delta,\xi}(\omega) = \frac{2\delta}{(1+\omega\cdot\xi)+\delta^2 (1-\omega\cdot\xi)} \ \left(\omega- (\omega\cdot\xi) \ \xi \right) 
+ \frac{(1+\omega\cdot\xi)-\delta^2 (1-\omega\cdot\xi)}{(1+\omega\cdot\xi)+\delta^2 (1-\omega\cdot\xi)} \ \xi \,.
\end{equation*}

\begin{lemma}
 \label{com}
Let $f\in L^1(\Sph^N)$ with $\int_{\Sph^N} f(\omega) \,d\omega\neq 0$. Then there is a transformation $\gamma_{\delta,\xi}$ of $\Sph^N$ such that
$$
\int_{\Sph^N} \gamma_{\delta,\xi}(\omega) f(\omega) \,d\omega = 0 \,.
$$
\end{lemma}

\begin{proof}
We may assume that $f\in L^1(\Sph^N)$ is normalized by $\int_{\Sph^N} f(\omega) \,d\omega=1$. We shall show that the $\R^{N+1}$-valued function
$$
F(r\xi) := \int_{\Sph^N} \gamma_{1-r,\xi}(\omega) f(\omega) \,d\omega\,,
\qquad 0< r< 1\,,\ \xi\in\Sph^N \,, 
$$
has a zero. First, note that because of $\gamma_{1,\xi}(\omega)= \omega$ for all $\xi$ and all $\omega$, the limit of $F(r\xi)$ as $r\to 0$ is independent of $\xi$. In other words, $F$ is a continuous function on the open unit ball of $\R^{N+1}$. In order to understand its boundary behavior, one easily checks that for any $\omega\neq-\xi$ one has $\lim_{\delta\to 0} \gamma_{\delta,\xi}(\omega)= \xi$, and that  this convergence is uniform on $\{(\omega,\xi) \in \Sph^N\times\Sph^N :\ 1+\omega\cdot\xi \geq \epsilon\}$ for any $\epsilon>0$. This implies that
\begin{equation*}
\lim_{r\to 1} F(r\xi) = \xi
\qquad\text{uniformly in}\ \xi\,.
\end{equation*}
Hence $F$ is a continuous function on the \emph{closed} unit ball, which is the identity on the boundary. The assertion is now a consequence of Brouwer's fixed point theorem.
\end{proof}

In the proof of Theorem \ref{mainjl} we use Lemma \ref{com} with $f=|u|^q$. Then the new function $\tilde u(\omega) = |J_{\gamma^{-1}}(\omega)|^{1/q} u(\gamma^{-1}(\omega))$, with $\gamma=\gamma_{\delta,\xi}$ of Lemma \ref{com}, satisfies the center of mass condition \eqref{eq:comjl}. Moreover, since rotations of the sphere, stereographic projection $\mathcal S$ and the dilations $\mathcal D_{\delta}$ leave the inequality invariant, $u$ can be replaced by $\tilde u$ in \eqref{eq:mainjl} without changing the values of each side.

In particular, if $U$ is an optimizer, our proof in Section~\ref{sec:jl} shows that the corresponding $\tilde U$ is a constant, which means that the original $U$ is a constant times $|J_\gamma|^{1/q}$. It is now a matter of computation, using the explicit form of $\gamma_{\delta,\xi}$, to verify that all such functions have the form of \eqref{eq:optjl}.

Conversely, let us verify that all the functions given in \eqref{eq:optjl} are optimizers. By the rotation invariance of inequality \eqref{eq:mainjl}, we can restrict our attention to the case $\xi= (0,\ldots,0,r)$ with $0<r<1$. These functions correspond via stereographic projection, \eqref{eq:corres}, to dilations of a constant times the function $H$ in \eqref{eq:opt}. Because of the dilation invariance of inequality \eqref{eq:sobjl} and because of the fact that we already know that $H$, which corresponds to the constant on the sphere, is an optimizer, we conclude that any function of the form \eqref{eq:optjl} is an optimizer.

We have discussed the derivative (Sobolev) version of the $\lambda=N-2$ case of \eqref{eq:mainsph}. Exactly the same considerations show the invariance of the fractional integral for all $0<\lambda<N$.


\subsection*{Acknowledgement} We thank Richard Bamler for valuable help with Appendix \ref{sec:com}.


\bibliographystyle{amsalpha}

\end{document}